\documentclass[12pt]{article}%
\usepackage{sw20bams}
\usepackage{amsmath}
\usepackage{amsfonts}
\usepackage{amssymb}
\usepackage{graphicx}%
\providecommand{\U}[1]{\protect\rule{.1in}{.1in}}
\begin{document}

\title{Pins, Stakes, Anchors and Gaussian Triangles}
\author{Steven Finch}
\date{October 24, 2014}
\maketitle

\begin{abstract}
Among the topics we discuss are certain joint densities (for sides and for
angles), acuteness probabilities and bivariate Rice moments.

\end{abstract}

\footnotetext{Copyright \copyright \ 2014 by Steven R. Finch. All rights
reserved.}A planar triangle is \textbf{pinned Gaussian} if one vertex is fixed
at the origin and the other two have $x$, $y$ coordinates that are independent
normally distributed random variables with mean $0$ and variance $1$. Let $a$,
$b$ denote the sides opposite the random vertices; let $c$ denote the side
opposite $(0,0)$. A\ Jacobian determinant calculation similar to
\cite{Fi1-PinStk} yields
\[
\left\{
\begin{array}
[c]{l}%
\dfrac{2}{\pi}\dfrac{x\,y\,z}{\sqrt{(x+y+z)(-x+y+z)(x-y+z)(x+y-z)}}\exp\left(
-\dfrac{1}{2}\left(  x^{2}+y^{2}\right)  \right) \\
\;\;\;\;\;\;\;\;\;\;\;\;\;\;\;\text{if }|x-y|<z<x+y,\\
0\;\;\;\;\;\;\;\;\;\;\;\;\;\text{otherwise}%
\end{array}
\right.
\]
as the trivariate density $f(x,y,z)$ for $a$, $b$, $c$. The condition
$|x-y|<z<x+y$ is equivalent to $|x-z|<y<x+z$ and to $|y-z|<x<y+z$ via the Law
of Cosines. As a consequence, the univariate density for $a$ (or $b$) is
\[%
\begin{array}
[c]{ccc}%
x\exp\left(  -\dfrac{1}{2}x^{2}\right)  , &  & x>0
\end{array}
\]
and the univariate density for $c$ is
\[%
\begin{array}
[c]{ccc}%
\dfrac{x}{2}\exp\left(  -\dfrac{1}{4}x^{2}\right)  , &  & x>0
\end{array}
\]
which are Rayleigh with means $\sqrt{\pi/2}$, $\sqrt{\pi}$ and mean squares
$2$, $4$ respectively. Side $c$ is distributed the same as that for
\textquotedblleft pure\textquotedblright\ Gaussian triangles.\ Also
$\operatorname*{E}(a\,b)=\pi/2$ and $\operatorname*{E}(a\,c)=\sqrt{2\pi}$.

From this, another calculation like \cite{Fi1-PinStk} gives
\[
\left\{
\begin{array}
[c]{lll}%
\dfrac{2}{\pi}\dfrac{\sin(x)\sin(y)\sin(x+y)}{\left(  \sin(x)^{2}+\sin
(y)^{2}\right)  ^{2}} &  & \text{if }0<x<\pi\text{, }0<y<\pi\text{ and
}x+y<\pi,\\
0 &  & \text{otherwise}%
\end{array}
\right.
\]
as the bivariate density for angles $\alpha$, $\beta$ opposite sides $a$, $b$.
It is clear that $\gamma=\pi-\alpha-\beta$ is Uniform[$0$,$\pi$] by
construction. \ Let $\wp$ denote the probability that a pinned Gaussian
triangle is obtuse. The univariate density for $\alpha$ (or $\beta$) is
obtained via integration:
\begin{align*}
&  \ \ \ \frac{2}{\pi}%
%TCIMACRO{\dint \limits_{0}^{\pi-x}}%
%BeginExpansion
{\displaystyle\int\limits_{0}^{\pi-x}}
%EndExpansion
\dfrac{\sin(x)\sin(y)\sin(x+y)}{\left(  \sin(x)^{2}+\sin(y)^{2}\right)  ^{2}%
}dy\\
\  &  =\tfrac{2}{\pi}%
%TCIMACRO{\dint \limits_{0}^{\pi-x}}%
%BeginExpansion
{\displaystyle\int\limits_{0}^{\pi-x}}
%EndExpansion
\tfrac{\cos(x)\sin(x)}{2\left(  2-\cos(x)^{2}\right)  \left(  \sin(x)^{2}%
+\sin(y)^{2}\right)  }dy\ \\
&  \ \ \ +\tfrac{2}{\pi}%
%TCIMACRO{\dint \limits_{0}^{\pi-x}}%
%BeginExpansion
{\displaystyle\int\limits_{0}^{\pi-x}}
%EndExpansion
\left(  \tfrac{\sin(x)\sin(y)\sin(x+y)}{\left(  \sin(x)^{2}+\sin
(y)^{2}\right)  ^{2}}-\tfrac{\cos(x)\sin(x)}{2\left(  2-\cos(x)^{2}\right)
\left(  \sin(x)^{2}+\sin(y)^{2}\right)  }\right)  dy\\
\  &  =\frac{1}{\pi}\frac{\cos(x)}{\left(  2-\cos(x)^{2}\right)  ^{3/2}%
}\left(  \frac{\pi}{2}+\arcsin\left(  \frac{\cos(x)}{\sqrt{2}}\right)
\right)  +\frac{1}{\pi}\frac{1}{2-\cos(x)^{2}}.
\end{align*}
Call this latter expression $g(x)$. Several results follow:%
\[%
\begin{array}
[c]{ccc}%
\operatorname*{E}(\alpha)=\dfrac{\pi}{4}, &  & \operatorname*{E}%
(\gamma)=\dfrac{\pi}{2},
\end{array}
\]%
\[%
\begin{array}
[c]{ccc}%
\operatorname*{E}(\alpha^{2})=\dfrac{5}{48}\pi^{2}+\dfrac{1}{4}\ln\left(
2\right)  ^{2}, &  & \operatorname*{E}(\gamma^{2})=\dfrac{\pi^{2}}{3},
\end{array}
\]%
\[%
\begin{array}
[c]{ccc}%
\operatorname*{E}(\alpha\,\beta)=\dfrac{1}{16}\pi^{2}-\dfrac{1}{4}\ln\left(
2\right)  ^{2}, &  & \operatorname*{E}(\alpha\,\gamma)=\dfrac{\pi^{2}}{12}.
\end{array}
\]
The cross-correlation coefficient
\[
\rho(\alpha,\beta)=\frac{\operatorname*{Cov}(\alpha,\beta)}{\sqrt
{\operatorname*{Var}(\alpha)\operatorname*{Var}(\beta)}}=\frac
{\operatorname*{E}(\alpha\,\beta)-\pi^{2}/16}{\operatorname*{E}(\alpha
^{2})-\pi^{2}/16}\approx-0.226
\]
indicates weak negative dependency; replacing $\beta$ by $\gamma$, such
association becomes stronger since $\rho(\alpha,\gamma)\approx-0.622$.
\ Finally,
\[
G(x)=%
%TCIMACRO{\dint \limits_{0}^{x}}%
%BeginExpansion
{\displaystyle\int\limits_{0}^{x}}
%EndExpansion
g(\xi)\,d\xi=\frac{1}{\pi}\frac{\sin(x)}{\left(  2-\cos(x)^{2}\right)  ^{1/2}%
}\left(  \frac{\pi}{2}+\arcsin\left(  \frac{\cos(x)}{\sqrt{2}}\right)
\right)  +\frac{1}{\pi}x
\]
which implies that
\begin{align*}
\wp &  =\operatorname*{P}(\alpha>\pi/2)+\operatorname*{P}(\beta>\pi
/2)+\operatorname*{P}(\gamma>\pi/2)\\
&  =2\left(  1-G(\pi/2)\right)  +1/2\\
&  =3/2-1/\sqrt{2}\approx0.793
\end{align*}
because a triangle can have at most one obtuse angle. \ See \cite{ES-PinStk,
Fi2-PinStk} for alternative approaches for computing $\wp$.

We conjecture that the bivariate density for $\alpha$, $\beta$ for pinned
Gaussian triangles in $n$-dimensional space is
\[
\left\{
\begin{array}
[c]{lll}%
C_{n}\dfrac{\sin(x)^{n-1}\sin(y)^{n-1}\sin(x+y)^{n-1}}{\left(  \sin
(x)^{2}+\sin(y)^{2}\right)  ^{n}} &  & \text{if }0<x<\pi\text{, }0<y<\pi\text{
and }0<x+y<\pi,\\
0 &  & \text{otherwise}%
\end{array}
\right.
\]
where $C_{n}=(n-1)2^{n-1}/\pi$. Values of $\wp$ for all $n$ are known
\cite{ES-PinStk} and are consistent with predictions based on our conjecture:
\[%
\begin{tabular}
[c]{|l|l|l|l|l|l|l|}\hline
$n$ & 3 & 4 & 5 & 6 & 7 & 8\\\hline
$\wp_{n}$ & $1-\frac{1}{\pi}$ & $\frac{3}{2}-\frac{5}{4\sqrt{2}}$ &
$1-\frac{4}{3\pi}$ & $\frac{3}{2}-\frac{43}{32\sqrt{2}}$ & $1-\frac{22}{15\pi}
$ & $\frac{3}{2}-\frac{177}{128\sqrt{2}}$\\\hline
\end{tabular}
\]
but the subtle dependency between angles remains open.

Henceforth let $c>0$ be constant. \ A planar triangle is \textbf{staked
Gaussian} if one vertex is fixed at $(0,0)$, another vertex is fixed at
$(c,0)$ and the third has $x$, $y$ coordinates that are independent
Normal($0$,$1$). The term \textit{stake} (as in \textquotedblleft staking a
tent\textquotedblright) is new in geometric probability, as far as is known.
\ For definiteness, define the vertices $A=(c,0)$, $B=(0,0)$ and $C=(u,v)$.
\ We shall enter into more details than previously, as the territory is
uncharted. \ Clearly%
\[%
\begin{array}
[c]{ccc}%
\tan(\alpha)=\dfrac{v}{c-u}, &  & \tan(\beta)=\dfrac{v}{u}.
\end{array}
\]
To compute the Jacobian of the transformation $(u,v)\mapsto(\alpha,\beta)$,
note that%
\[%
\begin{array}
[c]{lll}%
\sec(\alpha)^{2}\dfrac{\partial\alpha}{\partial u}=\dfrac{\partial}{\partial
u}\dfrac{v}{c-u}=\dfrac{v}{(c-u)^{2}}, &  & \sec(\alpha)^{2}\dfrac
{\partial\alpha}{\partial v}=\dfrac{\partial}{\partial v}\dfrac{v}{c-u}%
=\dfrac{1}{c-u},\\
\sec(\beta)^{2}\dfrac{\partial\beta}{\partial u}=\dfrac{\partial}{\partial
u}\dfrac{v}{u}=-\dfrac{v}{u^{2}}, &  & \sec(\beta)^{2}\dfrac{\partial\beta
}{\partial v}=\dfrac{\partial}{\partial v}\dfrac{v}{u}=\dfrac{1}{u}%
\end{array}
\]
and%
\[%
\begin{array}
[c]{ccc}%
\cos(\alpha)^{2}=\dfrac{(c-u)^{2}}{(c-u)^{2}+v^{2}}, &  & \cos(\beta
)^{2}=\dfrac{u^{2}}{u^{2}+v^{2}}%
\end{array}
\]
hence the determinant is%
\[
|J|=\left\vert
\begin{array}
[c]{cc}%
\dfrac{v}{(c-u)^{2}+v^{2}} & \dfrac{c-u}{(c-u)^{2}+v^{2}}\\
\dfrac{-v}{u^{2}+v^{2}} & \dfrac{u}{u^{2}+v^{2}}%
\end{array}
\right\vert =\dfrac{c\,v}{\left(  u^{2}+v^{2}\right)  \left[  (c-u)^{2}%
+v^{2}\right]  }.
\]
Solving for $u$, $v$ in terms of $\alpha$, $\beta$, we obtain%
\[%
\begin{array}
[c]{ccc}%
u=\dfrac{c\tan(\alpha)}{\tan(\alpha)+\tan(\beta)}, &  & v=\dfrac{c\tan
(\alpha)\tan(\beta)}{\tan(\alpha)+\tan(\beta)}.
\end{array}
\]
Substituting these expressions into the standard bivariate normal density
\[
\frac{1}{2\pi}\exp\left[  -\frac{1}{2}\left\{  u^{2}+v^{2}\right\}  \right]
\]
and dividing by $|J|$ yields%
\[
\frac{c^{2}}{2\pi}\exp\left[  -\frac{c^{2}}{2}\frac{\sin(\alpha)^{2}}%
{\sin(\alpha+\beta)^{2}}\right]  \frac{\sin(\alpha)\sin(\beta)}{\sin
(\alpha+\beta)^{3}}.
\]
Multiplying by $2$ gives the correct normalization. \ We have not attempted to
find any univariate densities or moments here. \ For simplicity, set $c=1$.
\ The acuteness probability can be found numerically: \ \
\[
1-\wp=\frac{1}{\pi}%
%TCIMACRO{\dint \limits_{0}^{\frac{\pi}{2}}}%
%BeginExpansion
{\displaystyle\int\limits_{0}^{\frac{\pi}{2}}}
%EndExpansion%
%TCIMACRO{\dint \limits_{\frac{\pi}{2}-\alpha}^{\frac{\pi}{2}}}%
%BeginExpansion
{\displaystyle\int\limits_{\frac{\pi}{2}-\alpha}^{\frac{\pi}{2}}}
%EndExpansion
\exp\left[  -\frac{1}{2}\frac{\sin(\alpha)^{2}}{\sin(\alpha+\beta)^{2}%
}\right]  \frac{\sin(\alpha)\sin(\beta)}{\sin(\alpha+\beta)^{3}}%
d\beta\,d\alpha=0.23685...
\]
and symbolically:%
\[
1-\wp=\frac{1}{2}\left[  -1+\operatorname{erf}\left(  \frac{1}{\sqrt{2}%
}\right)  +\frac{I_{0}(1/4)}{\exp(1/4)}\right]
\]
where $\operatorname{erf}$ is the error function and $I_{m}$ is the
$m^{\text{th}}$ modified Bessel function of the first kind. \ The exact
expression comes not from evaluating the double integral in some unforeseen
manner, but via geometry as follows.

Shift the Gaussian mean to the left by $1/2$. \ Likewise, translate the
triangle $ABC$ to $A^{\prime}B^{\prime}C^{\prime}$, where $A^{\prime}%
=(1/2,0)$, $B^{\prime}=(-1/2,0)$ and $C^{\prime}=(u-1/2,v)$. The set of points
$C$ where $\alpha+\beta=\pi/2$ is satisfied corresponds precisely to the
points $C^{\prime}=(x,y)$ on the circle given by%
\[%
\begin{array}
[c]{ccccc}%
x^{2}+y^{2}=1/4 &  & \text{or} &  & r=1/2
\end{array}
\]
(the latter is in polar coordinates). \ Acute triangles $ABC$ correspond to
points $C^{\prime}$ outside the circle but inside the strip $|x|=1/2$. \ The
probability that $C^{\prime}$ falls inside the strip is \
\[
\frac{1}{2\pi}%
%TCIMACRO{\dint \limits_{-1/2}^{1/2}}%
%BeginExpansion
{\displaystyle\int\limits_{-1/2}^{1/2}}
%EndExpansion
\;\;%
%TCIMACRO{\dint \limits_{-\infty}^{\infty}}%
%BeginExpansion
{\displaystyle\int\limits_{-\infty}^{\infty}}
%EndExpansion
\exp\left[  -\frac{1}{2}\left\{  \left(  x+\frac{1}{2}\right)  ^{2}%
+y^{2}\right\}  \right]  dy\,dx=\frac{1}{2}\operatorname{erf}\left(  \frac
{1}{\sqrt{2}}\right)
\]
and the probability that $C^{\prime}$ falls inside the circle is%
\begin{align*}
&  \frac{1}{2\pi}%
%TCIMACRO{\dint \limits_{-1/2}^{1/2}}%
%BeginExpansion
{\displaystyle\int\limits_{-1/2}^{1/2}}
%EndExpansion
\;\;%
%TCIMACRO{\dint \limits_{-\sqrt{1/4-x^{2}}}^{\sqrt{1/4-x^{2}}}}%
%BeginExpansion
{\displaystyle\int\limits_{-\sqrt{1/4-x^{2}}}^{\sqrt{1/4-x^{2}}}}
%EndExpansion
\exp\left[  -\frac{1}{2}\left\{  \left(  x+\frac{1}{2}\right)  ^{2}%
+y^{2}\right\}  \right]  dy\,dx\\
&  =\frac{1}{2\pi}%
%TCIMACRO{\dint \limits_{0}^{2\pi}}%
%BeginExpansion
{\displaystyle\int\limits_{0}^{2\pi}}
%EndExpansion
\,\,%
%TCIMACRO{\dint \limits_{0}^{1/2}}%
%BeginExpansion
{\displaystyle\int\limits_{0}^{1/2}}
%EndExpansion
\exp\left[  -\frac{1}{2}\left\{  r^{2}+r\cos(\theta)+\frac{1}{4}\right\}
\right]  r\,dr\,d\theta\\
&  =\exp(-1/8)%
%TCIMACRO{\dint \limits_{0}^{1/2}}%
%BeginExpansion
{\displaystyle\int\limits_{0}^{1/2}}
%EndExpansion
\exp\left(  -\frac{1}{2}r^{2}\right)  I_{0}\left(  \frac{r}{2}\right)  r\,dr\\
&  =\exp(-1/8)%
%TCIMACRO{\dint \limits_{0}^{1/8}}%
%BeginExpansion
{\displaystyle\int\limits_{0}^{1/8}}
%EndExpansion
\exp\left(  -s\right)  I_{0}\left(  2\sqrt{\frac{s}{8}}\right)  ds\\
&  =1-\frac{1}{2}\left[  1+\exp(-1/4)I_{0}(1/4)\right]  =\frac{1}{2}\left[
1-\frac{I_{0}(1/4)}{\exp(1/4)}\right]
\end{align*}
using a formula for what is called $J(1/8,1/8)$ in \cite{Go-PinStk,
Lu-PinStk}. \ This completes the calculation of $1-\wp$. \ 

It is shown in the Appendix that%
\[
\left\{
\begin{array}
[c]{l}%
\dfrac{2}{\pi}\dfrac{x\,y}{\sqrt{(x+y+c)(-x+y+c)(x-y+c)(x+y-c)}}\exp\left(
-\dfrac{1}{2}x^{2}\right) \\
\;\;\;\;\;\;\;\;\;\;\;\;\;\;\;\text{if }|x-y|<c<x+y,\\
0\;\;\;\;\;\;\;\;\;\;\;\;\;\text{otherwise}%
\end{array}
\right.
\]
is the bivariate density for sides $a$, $b$. The asymmetry is unsurprising.
Side $a$ is distributed the same as that for pinned\ Gaussian
triangles:\ Rayleigh with mean $\sqrt{\pi/2}$ and mean square $2$. For
simplicity, set $c=1$. \ The univariate density for $b$ is
\[%
\begin{array}
[c]{ccc}%
x\exp\left(  -\dfrac{1}{2}\left(  x^{2}+1\right)  \right)  I_{0}\left(
x\right)  , &  & x>0
\end{array}
\]
which is Rice with mean%
\[
\frac{\sqrt{2\pi}}{4e^{1/4}}\left(  3I_{0}\left(  \frac{1}{4}\right)
+I_{1}\left(  \frac{1}{4}\right)  \right)  =1.5485724605511453806302363...
\]
and mean square $3$. Also,
\begin{align*}
\operatorname*{E}(a\,b)  &  =\dfrac{2}{\pi}%
%TCIMACRO{\dint \limits_{0}^{\infty}}%
%BeginExpansion
{\displaystyle\int\limits_{0}^{\infty}}
%EndExpansion
x^{2}(x+1)\exp\left(  -\dfrac{1}{2}x^{2}\right)  E\left(  \frac{2\sqrt{x}%
}{x+1}\right)  dx\\
&  =2.2627965282687383013183035...
\end{align*}
where
\[
E(\xi)=%
%TCIMACRO{\dint \limits_{0}^{\pi/2}}%
%BeginExpansion
{\displaystyle\int\limits_{0}^{\pi/2}}
%EndExpansion
\sqrt{1-\xi^{2}\sin(\theta)^{2}}\,d\theta=%
%TCIMACRO{\dint \limits_{0}^{1}}%
%BeginExpansion
{\displaystyle\int\limits_{0}^{1}}
%EndExpansion
\sqrt{\dfrac{1-\xi^{2}t^{2}}{1-t^{2}}}\,dt
\]
is the complete elliptic integral of the second kind. \ An expression for
$\operatorname*{E}(a\,b)$ in terms of Meijer $G$-function values is possible
\cite{Gs-PinStk}.

A planar triangle is \textbf{anchored Gaussian} if one vertex is fixed at
$(-c/2,0)$, another vertex is fixed at $(c/2,0)$ and the third has $x$, $y$
coordinates that are independent Normal($0$,$1$). The term \textit{anchoring}
(as in \textquotedblleft anchoring a ship\textquotedblright) is again new in
this context. \ For definiteness, define the vertices $A=(c/2,0)$,
$B=(-c/2,0)$ and $C=(u,v)$. \ As before, here are details:%
\[%
\begin{array}
[c]{ccc}%
\tan(\alpha)=\dfrac{v}{\frac{c}{2}-u}=\dfrac{2v}{c-2u}, &  & \tan
(\beta)=\dfrac{v}{\frac{c}{2}+u}=\dfrac{2v}{c+2u};
\end{array}
\]%
\[%
\begin{array}
[c]{lll}%
\sec(\alpha)^{2}\dfrac{\partial\alpha}{\partial u}=\dfrac{\partial}{\partial
u}\dfrac{2v}{c-2u}=\dfrac{4v}{(c-2u)^{2}}, &  & \sec(\alpha)^{2}%
\dfrac{\partial\alpha}{\partial v}=\dfrac{\partial}{\partial v}\dfrac
{2v}{c-2u}=\dfrac{2}{c-2u},\\
\sec(\beta)^{2}\dfrac{\partial\beta}{\partial u}=\dfrac{\partial}{\partial
u}\dfrac{2v}{c+2u}=-\dfrac{4v}{(c+2u)^{2}}, &  & \sec(\beta)^{2}%
\dfrac{\partial\beta}{\partial v}=\dfrac{\partial}{\partial v}\dfrac{2v}%
{c+2u}=\dfrac{2}{c+2u};
\end{array}
\]%
\[%
\begin{array}
[c]{ccc}%
\cos(\alpha)^{2}=\dfrac{(c-2u)^{2}}{(c-2u)^{2}+4v^{2}}, &  & \cos(\beta
)^{2}=\dfrac{(c+2u)^{2}}{(c+2u)^{2}+4v^{2}}%
\end{array}
\]
hence the Jacobian determinant is%
\[
|J|=\left\vert
\begin{array}
[c]{cc}%
\dfrac{4v}{(c-2u)^{2}+4v^{2}} & \dfrac{2(c-2u)}{(c-2u)^{2}+4v^{2}}\\
-\dfrac{4v}{(c+2u)^{2}+4v^{2}} & \dfrac{2(c+2u)}{(c+2u)^{2}+4v^{2}}%
\end{array}
\right\vert =\dfrac{16c\,v}{\left[  (c-2u)^{2}+4v^{2}\right]  \left[
(c+2u)^{2}+4v^{2}\right]  }.
\]
Solving for $u$, $v$ in terms of $\alpha$, $\beta$, we obtain%
\[%
\begin{array}
[c]{ccc}%
u=\dfrac{c}{2}\dfrac{\tan(\alpha)-\tan(\beta)}{\tan(\alpha)+\tan(\beta)}, &  &
v=c\dfrac{\tan(\alpha)\tan(\beta)}{\tan(\alpha)+\tan(\beta)}.
\end{array}
\]
Substituting these expressions into the standard bivariate normal density and
dividing by $|J|$ yields%
\[
\frac{c^{2}}{2\pi}\exp\left[  -\frac{c^{2}}{8}\frac{\sin(\alpha-\beta
)^{2}+4\sin(\alpha)^{2}\sin(\beta)^{2}}{\sin(\alpha+\beta)^{2}}\right]
\frac{\sin(\alpha)\sin(\beta)}{\sin(\alpha+\beta)^{3}}.
\]
Multiplying by $2$ gives the correct normalization. \ We have not attempted to
find any univariate densities or moments here. \ For simplicity, set $c=1$.
\ The acuteness probability can be found numerically: \ \
\[
1-\wp=\frac{1}{\pi}%
%TCIMACRO{\dint \limits_{0}^{\frac{\pi}{2}}}%
%BeginExpansion
{\displaystyle\int\limits_{0}^{\frac{\pi}{2}}}
%EndExpansion%
%TCIMACRO{\dint \limits_{\frac{\pi}{2}-\alpha}^{\frac{\pi}{2}}}%
%BeginExpansion
{\displaystyle\int\limits_{\frac{\pi}{2}-\alpha}^{\frac{\pi}{2}}}
%EndExpansion
\exp\left[  -\frac{1}{8}\frac{\sin(\alpha-\beta)^{2}+4\sin(\alpha)^{2}%
\sin(\beta)^{2}}{\sin(\alpha+\beta)^{2}}\right]  \frac{\sin(\alpha)\sin
(\beta)}{\sin(\alpha+\beta)^{3}}d\beta\,d\alpha=0.26542...
\]
and symbolically:%
\[
1-\wp=-1+\frac{1}{\exp(1/8)}+\operatorname{erf}\left(  \frac{1}{2\sqrt{2}%
}\right)  .
\]
This value is slightly larger for anchored triangles than for staked
triangles. \ Proving this formula is similar to before, with $x+1/2$ replaced
by $x$ in both of the double integrals, which simplifies matters.

It is shown in the Appendix that%
\[
\left\{
\begin{array}
[c]{l}%
\dfrac{2e^{c^{2}/8}}{\pi}\dfrac{x\,y}{\sqrt{(x+y+c)(-x+y+c)(x-y+c)(x+y-c)}%
}\exp\left(  -\dfrac{1}{4}\left(  x^{2}+y^{2}\right)  \right) \\
\;\;\;\;\;\;\;\;\;\;\;\;\;\;\;\text{if }|x-y|<c<x+y,\\
0\;\;\;\;\;\;\;\;\;\;\;\;\;\text{otherwise}%
\end{array}
\right.
\]
is the bivariate density for sides $a$, $b$. The symmetry is very helpful. For
simplicity, set $c=1$. \ The univariate density for $a$ (or $b$) is
\[%
\begin{array}
[c]{ccc}%
x\exp\left(  -\dfrac{1}{2}\left(  x^{2}+\dfrac{1}{4}\right)  \right)
I_{0}\left(  \dfrac{x}{2}\right)  , &  & x>0
\end{array}
\]
which is Rice with mean%
\[
\frac{\sqrt{2\pi}}{16e^{1/16}}\left(  9I_{0}\left(  \frac{1}{16}\right)
+I_{1}\left(  \frac{1}{16}\right)  \right)  =1.3304473406107031708025583...
\]
and mean square $9/4$. Also,
\begin{align*}
\operatorname*{E}(a\,b)  &  =\tfrac{1}{64}(I_{0}\left(  \tfrac{1}{16}\right)
K_{0}\left(  \tfrac{1}{16}\right)  +8I_{0}\left(  \tfrac{1}{16}\right)
K_{1}\left(  \tfrac{1}{16}\right)  -8I_{1}\left(  \tfrac{1}{16}\right)
K_{0}\left(  \tfrac{1}{16}\right)  +I_{1}\left(  \tfrac{1}{16}\right)
K_{1}\left(  \tfrac{1}{16}\right)  )\\
&  =2.0303939030262620132069685...
\end{align*}
where $K_{m}$ is the $m^{\text{th}}$ modified Bessel function of the second kind.

The bivariate density for $\alpha$, $\beta$ / $a$, $b$ for staked/anchored
Gaussian triangles in $n$-dimensional space is also of interest, as well as
acuteness probabilities and cross-covariances.

In closing, let us return to \textquotedblleft pure\textquotedblright%
\ Gaussian triangles in $n$-dimensional space \cite{Fi1-PinStk, ES-PinStk}. We
conjecture that the bivariate density for $\alpha$, $\beta$ is \ \
\[
\left\{
\begin{array}
[c]{lll}%
\tilde{C}_{n}\dfrac{\sin(x)^{n-1}\sin(y)^{n-1}\sin(x+y)^{n-1}}{\left(
\sin(x)^{2}+\sin(y)^{2}+\sin(x+y)^{2}\right)  ^{n}} &  & \text{if }%
0<x<\pi\text{, }0<y<\pi\text{ and }x+y<\pi,\\
0 &  & \text{otherwise}%
\end{array}
\right.
\]
where $\tilde{C}_{n}=(n-1)2^{n-1}3^{n/2}/\pi$. \ Values of $\tilde{\wp}$ for
all $n$ are known \cite{ES-PinStk} and are consistent with predictions based
on our conjecture:
\[%
\begin{tabular}
[c]{|l|l|l|l|l|l|l|l|}\hline
$n$ & 2 & 3 & 4 & 5 & 6 & 7 & 8\\\hline
$\tilde{\wp}_{n}$ & $\frac{3}{4}$ & $1-\frac{3\sqrt{3}}{4\pi}$ & $\frac
{17}{32}$ & $1-\frac{9\sqrt{3}}{8\pi}$ & $\frac{353}{512}$ & $1-\frac
{27\sqrt{3}}{20\pi}$ & $\frac{867}{4096}$\\\hline
\end{tabular}
\]
Surely someone else has examined this issue! \ More relevant information would
be appreciated.

\section{Related Work}

Let $(X,Z)$ be bivariate normally distributed random variables with mean $0$,
variance $1$ and cross-correlation coefficient $\rho$. \ Let\ $(Y,W)$ be
likewise and independent of $(X,Z)$. The density of%
\[
(A,B)=\left(  \sqrt{(X-1/2)^{2}+Y^{2}},\sqrt{(Z-1/2)^{2}+W^{2}}\right)
\]
is \cite{Ml-PinStk, Si-PinStk}%
\[
\Omega\,a\,b\exp\left(  -\frac{a^{2}+b^{2}}{2(1-\rho^{2})}\right)
%TCIMACRO{\dsum \limits_{k=0}^{\infty}}%
%BeginExpansion
{\displaystyle\sum\limits_{k=0}^{\infty}}
%EndExpansion
\varepsilon_{k}I_{k}\left(  \frac{a\,b\,\rho}{1-\rho^{2}}\right)  I_{k}\left(
\frac{a}{2(1+\rho)}\right)  I_{k}\left(  \frac{b}{2(1+\rho)}\right)
\]
where $\varepsilon_{0}=1$, $\varepsilon_{k}=2$ for $k>0$ and the normalizing
constant is%
\[
\Omega=\frac{1}{1-\rho^{2}}\exp\left(  -\frac{1}{4(1+\rho)}\right)  .
\]
This result concerns the joint distribution of correlated distances from a
point $(1/2,0)$, which is not what we truly seek. \ Using an argument in
\cite{Si-PinStk} applied instead to a more general theorem in \cite{MB-PinStk}%
, we evaluate that
\[
(\bar{A},\bar{B})=\left(  \sqrt{(X+1/2)^{2}+Y^{2}},\sqrt{(Z-1/2)^{2}+W^{2}%
}\right)
\]
has density%
\[
\bar{\Omega}\,a\,b\exp\left(  -\frac{a^{2}+b^{2}}{2(1-\rho^{2})}\right)
%TCIMACRO{\dsum \limits_{k=0}^{\infty}}%
%BeginExpansion
{\displaystyle\sum\limits_{k=0}^{\infty}}
%EndExpansion
(-1)^{k}\varepsilon_{k}I_{k}\left(  \frac{a\,b\,\rho}{1-\rho^{2}}\right)
I_{k}\left(  \frac{a}{2(1-\rho)}\right)  I_{k}\left(  \frac{b}{2(1-\rho
)}\right)
\]
with normalizing constant%
\[
\bar{\Omega}=\frac{1}{1-\rho^{2}}\exp\left(  -\frac{1}{4(1-\rho)}\right)  .
\]
The changes from former to latter are slight. \ Our conjecture would be that
this joint distribution approaches that for anchored triangle sides (with
$c=1$) as $\rho\rightarrow1^{-}$.\ \ In particular, $\operatorname*{E}(\bar
{A}\,\bar{B})$ would approach $2.03039...$ as $\rho$ increases and
$\operatorname*{E}(\bar{A}\,\bar{B})$ can be rewritten as a doubly infinite
series as in \cite{Kr-PinStk}, avoiding integration entirely. \ Numerical
convergence issues associated with the doubly infinite series for large $\rho
$, however, make us hesitate to commit further.

\section{Acknowledgement}

I am thankful to M. Larry Glasser \cite{Gs-PinStk} for discovering the Meijer
$G$-function expression
\[
\operatorname*{E}(a\,b)=\frac{1}{4}\left[  G_{2,3}^{2,1}\left(  \frac{1}%
{2}\left\vert
\begin{array}
[c]{c}%
-\frac{3}{2},-\frac{1}{2}\\
-2,0,-2
\end{array}
\right.  \right)  +G_{2,3}^{2,1}\left(  \frac{1}{2}\left\vert
\begin{array}
[c]{c}%
-\frac{1}{2},-\frac{1}{2}\\
-1,0,-1
\end{array}
\right.  \right)  +G_{2,3}^{2,1}\left(  \frac{1}{2}\left\vert
\begin{array}
[c]{c}%
-\frac{1}{2},-\frac{1}{2}\\
0,0,-2
\end{array}
\right.  \right)  \right]
\]
corresponding to staked triangle sides. Much more relevant material can be
found at \cite{Fi3-PinStk}, including experimental computer runs that aided
theoretical discussion here.

\section{Appendix}

Let $\Delta=(a+b+c)(-a+b+c)(a-b+c)(a+b-c)$. \ The natural transformation
$(\alpha,\beta,c)\mapsto$ $(a,b,c)$ appearing in \cite{Fi1-PinStk} has
Jacobian determinant $a\,b$. Using the identities
\[%
\begin{array}
[c]{ccccc}%
\dfrac{a}{c}=\dfrac{\sin(\alpha)}{\sin(\alpha+\beta)}, &  & \dfrac{b}%
{c}=\dfrac{\sin(\beta)}{\sin(\alpha+\beta)}, &  & \dfrac{\sqrt{\Delta}}%
{2c^{2}}=\dfrac{\sin(\alpha)\sin(\beta)}{\sin(\alpha+\beta)}%
\end{array}
\]
we have%
\[
\sin(\alpha+\beta)=\frac{c}{a}\frac{c}{b}\dfrac{\sqrt{\Delta}}{2c^{2}}%
=\dfrac{\sqrt{\Delta}}{2a\,b}
\]
thus the bivariate staked angle density can be rewritten as%
\begin{align*}
&  \frac{c^{2}}{\pi}\exp\left[  -\frac{c^{2}}{2}\frac{\sin(\alpha)^{2}}%
{\sin(\alpha+\beta)^{2}}\right]  \frac{\sin(\alpha)\sin(\beta)}{\sin
(\alpha+\beta)^{3}}\frac{1}{a\,b}\\
&  =\frac{c^{2}}{\pi}\exp\left[  -\frac{c^{2}}{2}\frac{a^{2}}{c^{2}}\right]
\frac{a}{c}\frac{b}{c}\frac{2a\,b}{\sqrt{\Delta}}\frac{1}{a\,b}\\
&  =\dfrac{2}{\pi}\dfrac{a\,b}{\sqrt{\Delta}}\exp\left(  -\dfrac{1}{2}%
a^{2}\right)  .
\end{align*}
Also, we have%
\begin{align*}
\dfrac{\sin(\alpha-\beta)}{\sin(\alpha+\beta)}  &  =\dfrac{\sin(\alpha
)\cos(\beta)-\cos(\alpha)\sin(\beta)}{\sin(\alpha+\beta)}\\
&  =\frac{a}{c}\,\frac{-b^{2}+a^{2}+c^{2}}{2a\,c}-\frac{b}{c}\,\frac
{-a^{2}+b^{2}+c^{2}}{2b\,c}\\
&  =\frac{-b^{2}+a^{2}+c^{2}}{2c^{2}}-\frac{-a^{2}+b^{2}+c^{2}}{2c^{2}}\\
&  =\frac{a^{2}-b^{2}}{c^{2}}%
\end{align*}
thus the bivariate anchored angle density can be rewritten as%
\begin{align*}
&  \frac{c^{2}}{\pi}\exp\left[  -\frac{c^{2}}{8}\frac{\sin(\alpha-\beta
)^{2}+4\sin(\alpha)^{2}\sin(\beta)^{2}}{\sin(\alpha+\beta)^{2}}\right]
\frac{\sin(\alpha)\sin(\beta)}{\sin(\alpha+\beta)^{3}}\frac{1}{a\,b}\\
&  =\frac{c^{2}}{\pi}\exp\left[  -\frac{c^{2}}{8}\frac{\left(  a^{2}%
-b^{2}\right)  ^{2}+\Delta}{c^{4}}\right]  \frac{a}{c}\frac{b}{c}\frac
{2a\,b}{\sqrt{\Delta}}\frac{1}{a\,b}\\
&  =\frac{1}{\pi}\frac{2a\,b}{\sqrt{\Delta}}\exp\left[  -\frac{c^{2}}{8}%
\frac{\left(  a^{2}-b^{2}\right)  ^{2}+(a+b+c)(-a+b+c)(a-b+c)(a+b-c)}{c^{4}%
}\right] \\
&  =\frac{1}{\pi}\frac{2a\,b}{\sqrt{\Delta}}\exp\left[  -\frac{c^{2}}{8}%
\frac{c^{2}}{c^{4}}\left(  2a^{2}+2b^{2}-c^{2}\right)  \right] \\
&  =\frac{e^{c^{2}/8}}{\pi}\frac{2a\,b}{\sqrt{\Delta}}\exp\left[  -\frac{1}%
{4}\left(  a^{2}+b^{2}\right)  \right]  .
\end{align*}
Let $c=1$ from now on and let $\mathcal{S}$ denote the infinite strip
$\{(x,y):|x-y|<1<x+y\}$. \ The asymmetry within \
\[
\operatorname*{E}(a\,b)=\dfrac{2}{\pi}%
%TCIMACRO{\diint \limits_{\mathcal{S}}}%
%BeginExpansion
{\displaystyle\iint\limits_{\mathcal{S}}}
%EndExpansion
\dfrac{x^{2}y^{2}}{\sqrt{(x+y+1)(-x+y+1)(x-y+1)(x+y-1)}}\exp\left(  -\dfrac
{1}{2}x^{2}\right)  dy\,dx
\]
for staked triangles is actually advantageous since%
\[%
%TCIMACRO{\dint \limits_{\left\vert x-1\right\vert }^{x+1}}%
%BeginExpansion
{\displaystyle\int\limits_{\left\vert x-1\right\vert }^{x+1}}
%EndExpansion
\dfrac{y^{2}}{\sqrt{(x+y+1)(-x+y+1)(x-y+1)(x+y-1)}}dy=(x+1)E\left(
\frac{2\sqrt{x}}{x+1}\right)
\]
therefore a single (but difficult) integral remains. The symmetry within \
\[
\dfrac{2e^{1/8}}{\pi}%
%TCIMACRO{\diint \limits_{\mathcal{S}}}%
%BeginExpansion
{\displaystyle\iint\limits_{\mathcal{S}}}
%EndExpansion
\dfrac{x^{2}y^{2}}{\sqrt{(x+y+1)(-x+y+1)(x-y+1)(x+y-1)}}\exp\left(  -\dfrac
{1}{4}\left(  x^{2}+y^{2}\right)  \right)  dy\,dx
\]
for anchored triangles seems problematic at first, until we set%
\[%
\begin{array}
[c]{ccccccc}%
u=x+y, &  & v=-x+y; & \text{equivalently,} & 2\,x=u-v, &  & 2\,y=u+v
\end{array}
\]
and obtain%
\[%
\begin{array}
[c]{ccc}%
2\left(  x^{2}+y^{2}\right)  =u^{2}+v^{2}, &  & 4\,x\,y=u^{2}-v^{2}%
\end{array}
\]
therefore%
\begin{align*}
\operatorname*{E}(a\,b)  &  =\frac{1}{16}\dfrac{2e^{1/8}}{\pi}%
%TCIMACRO{\dint \limits_{1}^{\infty}}%
%BeginExpansion
{\displaystyle\int\limits_{1}^{\infty}}
%EndExpansion%
%TCIMACRO{\dint \limits_{-1}^{1}}%
%BeginExpansion
{\displaystyle\int\limits_{-1}^{1}}
%EndExpansion
\dfrac{\left(  u^{2}-v^{2}\right)  ^{2}}{\sqrt{\left(  u^{2}-1\right)  \left(
1-v^{2}\right)  }}\exp\left(  -\dfrac{1}{8}\left(  u^{2}+v^{2}\right)
\right)  \frac{1}{2}dv\,du\\
&  =\dfrac{e^{1/8}}{16\pi}%
%TCIMACRO{\dint \limits_{1}^{\infty}}%
%BeginExpansion
{\displaystyle\int\limits_{1}^{\infty}}
%EndExpansion
\dfrac{u^{4}}{\sqrt{u^{2}-1}}\exp\left(  -\dfrac{1}{8}u^{2}\right)  du\cdot%
%TCIMACRO{\dint \limits_{-1}^{1}}%
%BeginExpansion
{\displaystyle\int\limits_{-1}^{1}}
%EndExpansion
\dfrac{1}{\sqrt{1-v^{2}}}\exp\left(  -\dfrac{1}{8}v^{2}\right)  dv-\\
&  \dfrac{e^{1/8}}{8\pi}%
%TCIMACRO{\dint \limits_{1}^{\infty}}%
%BeginExpansion
{\displaystyle\int\limits_{1}^{\infty}}
%EndExpansion
\dfrac{u^{2}}{\sqrt{u^{2}-1}}\exp\left(  -\dfrac{1}{8}u^{2}\right)  du\cdot%
%TCIMACRO{\dint \limits_{-1}^{1}}%
%BeginExpansion
{\displaystyle\int\limits_{-1}^{1}}
%EndExpansion
\dfrac{v^{2}}{\sqrt{1-v^{2}}}\exp\left(  -\dfrac{1}{8}v^{2}\right)  dv+\\
&  \dfrac{e^{1/8}}{16\pi}%
%TCIMACRO{\dint \limits_{1}^{\infty}}%
%BeginExpansion
{\displaystyle\int\limits_{1}^{\infty}}
%EndExpansion
\dfrac{1}{\sqrt{u^{2}-1}}\exp\left(  -\dfrac{1}{8}u^{2}\right)  du\cdot%
%TCIMACRO{\dint \limits_{-1}^{1}}%
%BeginExpansion
{\displaystyle\int\limits_{-1}^{1}}
%EndExpansion
\dfrac{v^{4}}{\sqrt{1-v^{2}}}\exp\left(  -\dfrac{1}{8}v^{2}\right)  dv\\
&  =\frac{1}{64}\left[  K_{0}\left(  \dfrac{1}{16}\right)  +9K_{1}\left(
\dfrac{1}{16}\right)  \right]  I_{0}\left(  \dfrac{1}{16}\right)  -\\
&  \frac{1}{64}\left[  K_{0}\left(  \dfrac{1}{16}\right)  +K_{1}\left(
\dfrac{1}{16}\right)  \right]  \left[  I_{0}\left(  \dfrac{1}{16}\right)
-I_{1}\left(  \dfrac{1}{16}\right)  \right]  +\\
&  \frac{1}{64}K_{0}\left(  \dfrac{1}{16}\right)  \left[  I_{0}\left(
\dfrac{1}{16}\right)  -9I_{1}\left(  \dfrac{1}{16}\right)  \right]  ,
\end{align*}
as was to be shown.


\begin{thebibliography}{99}                                                                                               %


\bibitem {Fi1-PinStk}S. R. Finch, Random triangles, unpublished note (2010),
http://www.people.fas.harvard.edu/\symbol{126}sfinch/.

\bibitem {ES-PinStk}B. Eisenberg and R. Sullivan, Random triangles in $n$
dimensions, \textit{Amer. Math. Monthly} 103 (1996) 308--318; MR1383668 (96m:60025).

\bibitem {Fi2-PinStk}S. R. Finch, Random Gaussian tetrahedra, arXiv:1005.1033.

\bibitem {Go-PinStk}S. Goldstein, On the mathematics of exchange processes in
fixed columns. I, Mathematical solutions and asymptotic expansions,
\textit{Proc. Royal Soc. London. Ser. A} 219 (1953) 151--171; MR0058824 (15,429g).

\bibitem {Lu-PinStk}Y. L. Luke, \textit{Integrals of Bessel Functions},
McGraw-Hill, 1962, pp. 271--283; MR0141801 (25 \#5198).

\bibitem {Gs-PinStk}M. L. Glasser, private correspondence (2014).

\bibitem {Ml-PinStk}K. S. Miller, \textit{Multivariate Distributions}, Robert
E. Krieger Publishing Co., 1975, pp. 32--35; MR0381136 (52 \#2033).

\bibitem {Si-PinStk}M. K. Simon, Comments on \textquotedblleft Infinite-series
representations associated with the bivariate Rician distribution and their
applications\textquotedblright, \textit{IEEE Trans. Communications} 54 (2006)\ 1511--1512.

\bibitem {MB-PinStk}K. S. Miller, R. I. Bernstein and L. E. Blumenson,
Generalized Rayleigh processes, \textit{Quart. Appl. Math.} 16 (1958)
137--145; addendum 20 (1963) 395; MR0094862 (20 \#1371).

\bibitem {Kr-PinStk}M. Krishnan, The noncentral bivariate chi distribution,
\textit{SIAM Review} 9 (1967) 708--714; MR0224191 (36 \#7235).

\bibitem {Fi3-PinStk}S. R. Finch, Simulations in R\ involving triangles and
tetrahedra, http://www.people.fas.harvard.edu/\symbol{126}sfinch/csolve/rsimul.html.
\end{thebibliography}
\end{document}